\documentclass[11pt]{article}

\usepackage{amsfonts}
\usepackage{amsthm}
\usepackage{amsmath}
\usepackage{epsfig}
\usepackage{color}
\usepackage{hyperref}
\usepackage[all,tips]{xy}

\DeclareMathOperator{\vol}{vol}
\DeclareMathOperator{\PSL}{PSL} \DeclareMathOperator{\SL}{SL}
\DeclareMathOperator{\Ram}{Ram} \DeclareMathOperator{\Stab}{Stab}
\DeclareMathOperator{\PGL}{PGL}

\def\qed{ $\sqcup\!\!\!\!\sqcap$}

\def\nr{\rm{nr}}

\newcommand{\innp}[1]{\left< #1 \right>}
\newcommand{\abs}[1]{\left\vert#1\right\vert}
\newcommand{\set}[1]{\left\{#1\right\}}

\newcommand{\pr}[1]{\left( #1 \right) }

\newtheorem{theorem}{Theorem}[section]
\newtheorem{lemma}[theorem]{Lemma}

\newtheorem{proposition}[theorem]{Proposition}
\newtheorem{question}{Question}

\title{\textbf{The genus spectrum of a hyperbolic 3--manifold} }
\author{D.~B. McReynolds\thanks{This work was supported in part by an NSF
postdoctoral fellowship} ~and A.~W. Reid\thanks{This work was
supported in part by the NSF}}

\begin{document}
\maketitle

%---------------------------------------------------------------
%---------------------------------------------------------------
\section{Introduction}

Let $M$ be a closed orientable Riemannian manifold $M$. An invariant
of $M$ that has traditionally been of great importance in
understanding the geometry of $M$ is the geodesic length spectrum
(that is the set of lengths of closed geodesics counted with
multiplicities). For example, in the case when $M$ admits a metric
with all sectional curvatures negative, there is a strong
relationship between this spectrum and the eigenvalue spectrum of
the Laplace--Beltrami operator (\cite{DG},\cite{Gan});
the latter is well known
to determine the volume and dimension of $M$. The aim of this
article is to begin the development of more general geometric
spectra for Riemannian manifolds.

In more detail, the geodesic length spectrum encodes isometric
immersions of $S^1$ into $M$, and in this article we will study the
2--dimensional case; that is totally geodesic immersions of
orientable, finite type surfaces into Riemannian manifolds. Motivated by the case
of the length spectrum two natural (complementary) questions are:

\begin{question}
\label{questione1} How much of the geometry of $M$ is determined by
the totally geodesic surfaces immersed in $M$?
\end{question}

\begin{question}
\label{question2} Do there exist non-isometric Riemannian manifolds
with the same ``spectra'' of totally geodesic surfaces?
\end{question}

To state our results we describe the set up more carefully. We will
consider only the case when the manifold $M$ is a complete
orientable hyperbolic 3--manifold of finite volume. Let
$\Sigma_{g,n}$ denote the closed orientable surface of genus $g$
with $n$ punctures and let $M$ be as indicated above. For each
finite type hyperbolic surface $X$ in the moduli space of such a
surface $\Sigma_{g,n}$, set $\ell_X(M)$ to be the number of the free homotopy
classes of totally geodesic surfaces in $M$ with a representative
isometric to $X$. That $\ell_X(M)$ is finite follows from \cite{T}
Corollary 8.8.6. We define the {\em geometric genus spectrum} to be
the set of pairs
\[ {\cal GS}(M) = \{(X,\ell_X(M))~:~\ell_X(M) \ne 0\}. \]

We say that $M_1$ and $M_2$ are {\em geometrically isogenus} if
${\cal GS}(M_1) = {\cal GS}(M_2)$. Forgetting the multiplicities
$\ell_X(M)$ in $\mathcal{GS}(M)$ leads to the {\em geometric genus
set} given by
\[\textrm{GS}(M) = \{X~:~\ell_X(M) \ne 0\}. \]
We say that $M_1$ and $M_2$ are {\em geometrically genus equivalent} if $\textrm{GS}(M_1)=\textrm{GS}(M_2)$.

Our first result addresses the degree to which the geometric genus
set or spectrum governs the geometry of $M$.

\begin{theorem}
\label{main3}
Let $M_1={\bf H}^3/\Gamma_1$ and $M_2={\bf H}^3/\Gamma_2$ be arithmetic hyperbolic 3--manifolds. If $M_1$ and $M_2$ are
geometrically genus equivalent then either $\textrm{GS}(M_1)=\emptyset$,
or $M_1$ and $M_2$ are commensurable.
\end{theorem}

The first possibility does occur since, unlike the case of the length spectrum, the set
${\cal GS}(M)$ can be empty; indeed, \emph{most} finite volume
hyperbolic 3--manifolds do not contain an immersed totally geodesic
surface (see \cite{MR} Chapter 5.3). However, and perhaps most
interestingly from the point of view of this spectrum, if an
arithmetic hyperbolic 3--manifold contains one totally geodesic
surface, it contains infinitely many such surfaces (up to
commensurability). We do not know if this holds in general (see
Section \ref{Final} for a further discussion of this). Thus
arithmetic manifolds provide a good class to analyze in regard to geometric spectra. It is also
worth remarking that even in the case of the length spectrum, there
are no comparable results to Theorem \ref{main3} without an
arithmetic assumption---see \cite{Re}, \cite{CHLR}, and \cite{PR}.

Similarities with results for the length (and eigenvalue) spectrum
are also evident in our two other main results (cf. \cite{Sun},
\cite{LMNR06} and \cite{PR}).

\begin{theorem}
\label{main2} Let $M=\mathbf{H}^3/\Gamma$ be an arithmetic
hyperbolic $3$--manifold. Then there exists infinitely many pairs of
finite covers $(M_j,N_j)$ of $M$ such that $M_j,N_j$ are
geometrically isogenus and nonisometric.
\end{theorem}

\begin{theorem}
\label{main1} Let $M =\mathbf{H}^3/\Gamma$ be an arithmetic
hyperbolic $3$--manifold. Then there exists infinitely many pairs of
finite covers $(M_j,N_j)$ of $M$ such that $M_j,N_j$ are
geometrically genus equivalent and the sequence
$\set{\vol(M_j)/\vol(N_j)}_j$ is unbounded.
\end{theorem}

Though the method of proof for Theorems \ref{main2} and \ref{main1} is in the
same spirit as Sunada's method for producing iso-length and
isospectral manifolds \cite{Sun}, unlike Sunada's method it is not
purely algebraic (see the discussion in Section \ref{WithoutVolume}
for more on this comparison).

Theorems \ref{main2} and \ref{main1} also hold without an arithmetic
assumption. However, we have elected only to completely treat the case of
arithmetic hyperbolic 3--manifolds since, as mentioned above, it is
only in this case that we know that the geometric genus spectrum
being nontrivial implies that there are infinitely many totally
geodesic surfaces.  In addition the proofs in the arithmetic cases
are somewhat simpler. In Section \ref{Non-Arithmetic} we sketch the
modifications needed in the proofs for non-arithmetic manifolds.

%---------------------------------------------------------------
\subsection*{Acknowledgements}

Both authors wish to thank R.I.M.S. Kyoto for their hospitality in
December 2006 when this work began. Additionally, the first author
wishes to thank the California Institute of Technology as some of
this work was done during an extended visit to that institution, and the
second author wishes to thank the Institute for Advanced Study where this
work was completed.

%---------------------------------------------------------------
%---------------------------------------------------------------
\section{Preliminaries}

In this section we collect some preliminary material and notation
that we will use throughout. For convenience, we will often blur the
distinction between Kleinian groups as subgroups of $\SL(2,{\bf C})$
and $\PSL(2,{\bf C})$.\smallskip\smallskip

\noindent {\bf Notation:}~For a number field $K$, $R_K$ will denote
the ring of $K$--integers. If $L/K$ is an extension of number
fields, $P$ a prime ideal of $R_K$, and $\mathfrak{p}$ an ideal of
$R_L$ such that $\mathfrak{p}|P$, we then denote by
$\omega_P$ (resp. $\nu_{\mathfrak{p}}$) the places
associated to these primes and write $\nu_{\mathfrak{p}}|\omega_P$.
For a quadratic extension $K'/K$ of number fields and a $K$--prime
$\mathfrak{p}$, by the {\em splitting type} of $\mathfrak{p}$ in
$K'$, we mean whether $\mathfrak{p}$ is ramified, inert, or split in
$K'/K$. If $B$ is a quaternion algebra over $K$, we denote the set
of finite places of $K$ at which $B$ is ramified by $\Ram_f B$.

%---------------------------------------------------------------
\subsection{}\label{FieldsSec} 

For a prime $p\in {\bf N}$, $\mathbf{F}_q$
will denote the unique field of order
$q=p^j$ for each $j\geq 1$ and $\mathbf{F}_q(\sqrt{b})$ the unique quadratic extension
of $\mathbf{F}_q$. In a slight abuse of
notation, the undecorated $b$ above represents some element of
$\mathbf{F}_q$ that is not a square in $\mathbf{F}_q$, in essence
suppressing its dependence on the $p$ and $j$. As no specific
property of $b$ will be used, this abuse is mild. In what follows,
to avoid some complications, we will assume $p \ne 2$.

Associated to each field $\mathbf{F}_{q}$ is the finite group
$\PSL(2,\mathbf{F}_{q})$ which is well known to have order
$\frac{q(q^2-1)}{2}$. The standard \emph{Borel subgroup} of
$\PSL(2,\mathbf{F}_{q})$ is given by
\[ \mathbf{B}(\mathbf{F}_{q}) = \set{ \begin{pmatrix} \alpha &
\beta \\ 0 & \alpha^{-1} \end{pmatrix} ~:~ \alpha \in
(\mathbf{F}_{q})^\times/\set{\pm 1},~ \beta \in \mathbf{F}_{q} }\]
and has order $\frac{q(q-1)}{2}$.

In our proofs of Theorems \ref{main2} and \ref{main1}, we will need
elements of a prescribed order. This is accomplished with the
following well-known lemma.

\begin{lemma}\label{OrderLemma}
For each divisor $m$ of $(q\pm 1)$, there exists $\gamma_m$ in
$\SL(2,\mathbf{F}_q)$ with order $m$.
\end{lemma}

\noindent \textbf{Proof:} For divisors of $q-1$, simply note that we
have the inclusion of $(\mathbf{F}_{q})^\times$ into
$\SL(2,\mathbf{F}_{q})$ given by
\[ \alpha \mapsto \begin{pmatrix} \alpha & 0 \\ 0 & \alpha^{-1}
\end{pmatrix}. \]
The existence of $\gamma_m$ now follows from the fact that
$(\mathbf{F}_{q})^\times$ is cyclic.

For divisors of $q+1$, the ensuing argument produces the desired
elements. As right multiplication of
$(\mathbf{F}_{q}(\sqrt{b}))^\times$ on $\mathbf{F}_{q}(\sqrt{b})$ is
$\mathbf{F}_{q}$--linear, there exists an injective homomorphism
$(\mathbf{F}_{q}(\sqrt{b}))^\times$ into
$\textrm{GL}(2,\mathbf{F}_{q})$ given by selecting an
$\mathbf{F}_q$--basis for $\mathbf{F}_q(\sqrt{b})$. For any element
$\gamma$ in $(\mathbf{F}_q(\sqrt{b}))^\times$, $\det (\gamma^{q-1}) =
1$ since $\det \gamma \in \mathbf{F}_q^\times$ and determinants are multiplicative. Appealing again to
the fact that $(\mathbf{F}_q(\sqrt{b}))^\times$ is cyclic, we may
take $\gamma$ of order $q^2-1$ and thus produce an element
$\gamma^{s(q-1)}$ in $\SL(2,\mathbf{F}_q)$ of order $m$ where $ms =
q+1$. \qed

%---------------------------------------------------------------
\subsection{} 

We will require some control over the behavior of
primes in extensions. The first result of this type we
need is on the splitting of primes in quadratic extensions.

\begin{proposition}
\label{sameplittype} Let $K_1$ and $K_2$ be quadratic extensions of
a number field $K$. If all but a finite number of $K$--primes have
the same splitting type in $K_1$ and $K_2$, then
$K_1=K_2$.\end{proposition}

\noindent{\bf Proof:}~For $i=1,2$, let $S_i$ be the set of primes of
$K$ that split in $K_i$. Since $S_1=S_2$ apart from a finite number
of primes (in particular density 0), we deduce from \cite{Ja}
Chap. IV Corollary 5.5 that $K_1=K_2$ as required.\qed\smallskip\smallskip

Now let $k/L$ be a proper extension of number fields, $\omega$ a
finite place of $L$ and $\nu |\omega$ a place of $k$. For brevity we
simply denote the extension of residue class fields by
$\mathbf{F}_{q_\nu}/\mathbf{F}_{q_\omega}$.

\begin{lemma}\label{Ceb}
There exists a positive density set of finite places $\nu$ of $k$ such
that $\mathbf{F}_{q_\nu}/\mathbf{F}_{q_\omega}$ is an extension of
degree at least $2$. Moreover, we can arrange for both
$[\mathbf{F}_{q_\nu}:\mathbf{F}_{q_\omega}]$ and
$[\mathbf{F}_{q_\omega}:\mathbf{F}_{p_\nu}]$ to be fixed over this collection of places.
\end{lemma}

\noindent \textbf{Proof:} Let $k_{\textrm{gal}}$ denote the Galois
closure of $k$. As $k \ne L$, there exists a Galois automorphism
$\theta$ in $\textrm{Gal}(k_{\textrm{gal}}/L)$ such that $\theta$ is
trivial on $L$ and nontrivial on $k$. According to the Cebotarev Density Theorem, there exists
a positive density set of places $\nu_{\textrm{gal}}$ of
$k_{\textrm{gal}}$ such that the induced Galois automorphism
$\overline{\theta}$ on the Galois extension
$\mathbf{F}_{q_{\nu_{\textrm{gal}}}}$ of $\mathbf{F}_{q_\omega}$ has
order $\abs{\theta}$. In particular, $\overline{\theta}$ is trivial
on $\mathbf{F}_{q_\omega}$ but not $\mathbf{F}_{q_\nu}$. Thus
$\mathbf{F}_{q_\nu}/\mathbf{F}_{q_\omega}$ is a proper extension for
all the residue fields $\mathbf{F}_{q_\nu},\mathbf{F}_{q_\omega}$
associated to the places $\nu$ and $\omega$ where
$\nu_{\textrm{gal}} \mid \nu$ and $\nu_{\textrm{gal}} \mid \omega$.
For the final claim, let ${\bf F}_{p_\nu}$ denote the prime field
for $\mathbf{F}_{q_\nu}$ (and hence for $\mathbf{F}_{q_\omega}$ too). As the
possibilities for the degrees
$[\mathbf{F}_{q_\nu}:\mathbf{F}_{q_\omega}]$ and
$[\mathbf{F}_{q_\omega}:{\bf F}_{p_\nu}]$ range over finite sets as
we vary $\nu$, the result follows immediately from the Pigeon Hole
Principle and the fact that a finite union of density zero sets is
density zero. Specifically, at least one possibility for
$[\mathbf{F}_{q_\nu}:\mathbf{F}_{q_\omega}]$ and
$[\mathbf{F}_{q_\omega}:{\bf F}_{p_\nu}]$ with
$\mathbf{F}_{q_\nu}/\mathbf{F}_{q_\omega}$ proper occurs with
positive density.\qed

%---------------------------------------------------------------
\subsection{}\label{MaxFuchsianSect}

We recall the basic framework of non-elementary Fuchsian subgroups
of Kleinian groups. We will not restrict to the arithmetic setting at present,
as some of this material will be required for the general discussion in
Section \ref{Non-Arithmetic}.

Let $M={\bf H}^3/\Gamma$ be a finite volume hyperbolic 3--manifold
and $S$ a properly immersed orientable totally geodesic surface in
$M$. Then corresponding to $S$ is a non-elementary Fuchsian subgroup
of $F<\Gamma$ and a round circle ${\cal C} \subset \widehat{{\bf C}} = {\bf C} \cup \infty$
such that $F<\Stab({\cal C},\Gamma)$ where
$$\Stab({\cal C},\Gamma) = \{\gamma\in\Gamma : \gamma {\cal C} = {\cal C}~
\hbox{and}~\gamma~\hbox{preserves components of}~ \widehat{{\bf C}}\setminus {\cal C}\}.$$
We refer to the groups $\Stab({\cal C},\Gamma)$ as {\em
maximal Fuchsian subgroups} of $\Gamma$. Note that the geometric genus
spectrum of $M$ is determined by the maximal Fuchsian subgroups of $\Gamma$.

If $F$ is a non-elementary Fuchsian subgroup of a finite co-volume
Kleinian group $\Gamma$, then the invariant trace-field $L=kF$ is a
real number field and is a subfield of the invariant trace-field
$k=k\Gamma$. Furthermore the invariant algebra $AF/L$ is a
$L$--subalgebra of $A\Gamma$ and $AF\otimes_{L} k \cong A\Gamma$
as $k$--quaternion algebras.

Suppose that $A\Gamma$ is ramified at a finite place $\nu$ and
$\omega$ is a place of $kF$ such that $\nu |\omega$. Then $AF$ is
ramified at $\omega$ since
$$A\Gamma\otimes_k k_\nu \cong (AF\otimes_L k) \otimes_k k_\nu$$
is isomorphic to $(AF\otimes_L {L_\omega}) \otimes_{L_\omega}
k_\nu$ as $k_\nu$--quaternion algebras, and so $AF\otimes_L {L_\omega}$ must be a division algebra.
However, $AF$ can be ramified at other places (see for example the
arithmetic case below).  On the other hand, if $AF$ is unramified at
$\omega$ then it is clear that $A\Gamma$ is unramified at any
$k$--place $\nu$ that divides $\omega$.

Let $\nu$ be a place of $k$ at which $A\Gamma$ is unramified. Assume
that $\Gamma \subset A\Gamma^1$ (the group of elements of norm 1 in
$A\Gamma$) and assume additionally that under the induced injection
of $\Gamma$ into $\textrm{M}(2,k_\nu)$, the image is conjugate into
$\SL(2,R_\nu)$ where $R_\nu$ is the ring of $\nu$--adic integers in
$k_\nu$. Thus we will assume that $\Gamma$ injects into
$\SL(2,R_\nu)$. Note that this is true for all but a finite number of places
$\nu$.

Reducing modulo the maximal ideal $\pi R_\nu$ of $R_\nu$ and taking
the central quotient, gives a homomorphism $\phi_\nu\colon \Gamma
\rightarrow \PSL(2,\mathbf{F}_{q_\nu})$. For all but a finite number
of places, the Strong Approximation Theorem (see \cite{Nori87} or
\cite{Weis}) ensures that $\phi_\nu$ is surjective. For a place
$\omega$ of $L$ with $\nu | \omega$, the image of $F$ under
$\phi_\nu$ depends on whether the $AF$ is ramified at the place $\omega$ or
not. More specifically, we will prove the following.

\begin{lemma}\label{FuchsianImage}
In the notation established above, the image of $F$ under $\phi_\nu$ is conjugate (over an algebraic
closure of $\mathbf{F}_{q_\nu}$) to either
\begin{itemize}
\item[(a)] a subgroup of $\PSL(2,\mathbf{F}_{q_\omega})$, when $AF$ is not ramified at $\omega$, or
\item[(b)]
a subgroup of $\mathbf{B}(\mathbf{F}_{q_\omega}(\sqrt{b}))$ when $AF$ is ramified at $\omega$.
\end{itemize}
\end{lemma}

\noindent {\bf Proof:}~Suppose first that $AF$ is unramified at $\omega$. In this case, $AF \otimes_L L_\omega$
is isomorphic to $M(2,L_\omega)$, and it is
easy to see  from this that the image of
$F$ under $\phi_\nu$ is
isomorphic to a subgroup of $\PSL(2,\mathbf{F}_{q_\omega})$---indeed it will be conjugate over some finite
extension of $\mathbf{F}_{q_\omega}$ to a subgroup of $\PSL(2,\mathbf{F}_{q_\omega})$.

The case when $AF$ is ramified at $\omega$ requires some additional discussion. It will be helpful to
recall a certain representation of the division algebra $AF \otimes_L L_\omega$ (see
\cite{MR} Exercise 2.6(1)). Let $L_{\nr}$ denote the
unique unramified quadratic extension of $L_\omega$, and $'$ denote the non-trivial
Galois automorphism of $L_{\nr}/L$.  Then the unique division algebra over $L_\omega$ can be
represented as:

$$B_\omega=\{\begin{pmatrix} a & b \\ \pi_\omega b' &  a'\end{pmatrix} : a,b\in L\}.$$

Hence, $F$ is conjugate (over an algebraic closure of $L_\omega$
and hence $k_\nu$) to a subgroup of the elements of norm 1 in $B_\omega$.
Indeed, by assumption (see the discussion preceding the lemma), since $F$ is a subgroup
of
the unique maximal order of $AF \otimes_L L_\omega$ (see \cite{MR} Chap.6.4), it follows that
$F$ is actually conjugate to a subgroup of the norm 1 elements of
${\cal O}_\omega$, the unique maximal order in $B_\omega$. This discussion, together with
\cite{MR} Exercise 6.4(1), shows that $F$ is conjugate to a subgroup of

$$\{\begin{pmatrix} a & b \\ \pi_\omega b' &  a'\end{pmatrix} : a,b\in R\},$$
where $R$ is the valuation ring of $L_{\nr}$.

Now observe that the image of $F$ under the homomorphism $\phi_\nu$
is, up to conjugacy, described by its image under the natural reduction quotient of
${\cal O}_\omega$ by its unique maximal ideal (which is informally, the reduction
modulo $\pi_\omega$). This can be seen to have image in a group of upper triangular matrices with entries
in the unique quadratic extension of $\mathbf{F}_{q_\omega}$. This completes the
proof.\qed
%---------------------------------------------------------------
%---------------------------------------------------------------
%---------------------------------------------------------------
\section{Commensurability class rigidity: Theorem \ref{main3}}

Here we will prove Theorem \ref{main3}. This should be compared with
the results found in \cite{Re}, \cite{CHLR}, \cite{LSV}, and
\cite{PR} for the length and eigenvalue spectra.

%---------------------------------------------------------------
\subsection{}

We begin by specializing the discussion of Subsection
\ref{MaxFuchsianSect} to the case of arithmetic Kleinian groups
which contain Fuchsian subgroups. For a fuller account see \cite{MR}
Chapters 5.3 and 9.5.

Let $\Gamma$ be an arithmetic Kleinian group which contains a
non-elementary Fuchsian subgroup $F$.  It is known in this case that
if $\cal C$ is the invariant circle of $F$, then $\Stab({\cal
C},\Gamma)$ is a Fuchsian group of finite co-area. Furthermore, the
invariant trace-field $L$ of $F$ is a totally real field, $[k\Gamma
:L]=2$, and if $B$ is the invariant quaternion algebra of $F$ then
(as above), $B\otimes_{L} k\Gamma = A\Gamma$.  However, more can be
said in this case about precisely which algebras $B$ can arise. We
record the following (see \cite{MR} Theorem 9.5.5) which will be useful
in the sequel.

\begin{theorem}
\label{fuchsiansubgroups} Let $\Gamma$ be an arithmetic Kleinian
group with invariant trace-field $k$ and invariant quaternion
algebra $A/k$. Suppose $[k:L]=2$ where $L = k\cap{\bf R}$ and $B$ is
a quaternion algebra over $L$ ramified at all real places of $L$
except the identity. Then $A\cong B\otimes_L k$ if and only if
$\Ram_f A$ consists of $2r$ places (possibly zero)
$\{\nu_{{\mathfrak{p} }_1},\nu_{{\mathfrak{p} }_1'},\ldots
\nu_{{\mathfrak{p} }_r},\nu_{{\mathfrak{p} }_r'}\}$ where the
$k$--primes ${\mathfrak{p} }_i$ and ${\mathfrak{p}}_i'$ satisfy
${\mathfrak{p}}_i\cap R_L = {\mathfrak{p}}_i'\cap R_L = P_i$,
$\{\omega_{P_1},\ldots , \omega_{P_r}\} \subset \Ram_f B$ with
$\Ram_f B\setminus \{\omega_{P_1},\ldots ,\omega_{P_r}\}$ consisting
of primes in $R_L$ which are inert or ramified in
$k/L$.\end{theorem}

%---------------------------------------------------------------
\subsection{Proof of Theorem \ref{main3}}

Recall in the statement of Theorem \ref{main3}, we are assuming that we have a pair
of arithmetic hyperbolic 3--manifolds $M_1=\mathbf{H}^3/\Gamma_1$ and
$M_2=\mathbf{H}^3/\Gamma_2$ that are geometrically genus equivalent. Our goal
is to show that if $\textrm{GS}(M_1) \ne \emptyset$, then $M_1$ and $M_2$ are commensurable. For the proof, set $k_1, k_2$
and $A_1/k_1, A_2/k_2$ to be the invariant trace-fields and quaternion
algebras of $\Gamma_1,\Gamma_2$ respectively. With this notation, to show that $M_1$ and
$M_2$ are commensurable it suffices to show that $k_1 \cong k_2$ and
$A_1 \cong A_2$ as algebras over $\bf Q$ (see \cite{MR} Theorem
8.4.1).

Since $\Gamma_1$ contains a non-elementary Fuchsian subgroup, $k_1$
contains a totally real field $L$ with $[k_1:L]=2$. Furthermore,
since $M_1$ and $M_2$ are geometrically genus equivalent, it follows
from \S 3.1, that $L$ is the maximal
totally real field of $k_2$  and that  $[k_2:L]=2$ also holds.

Any non-elementary Fuchsian subgroup of $\Gamma_j$ defines an
associated quaternion algebra $B$ over $L$ that is an
$L$--subalgebra of $A_j$. Conversely, any quaternion algebra $B/L$
that is unramified at the identity, ramified at all other real
places of $L$, and for which $B \otimes_L k_j \cong A_j$, produces a
commensurability class of arithmetic Fuchsian subgroups of
$\Gamma_j$ (note that different embeddings of $B$ into $A_j$ can
provide non-conjugate subgroups of $\Gamma_j$). Since $M_1$ and
$M_2$ are geometrically genus equivalent, it follows that if $B$ is
a quaternion algebra defined over $L$ that is unramified at the
identity and ramified at all other real places of $L$, then:
\begin{equation}\label{1}
B \otimes_L k_1 \cong A_1~~\hbox{if and only if}~~B \otimes_L k_2
\cong A_2.
\end{equation}
In particular, it is easy to see that  if $k_1 \cong k_2$, then $A_1
\cong A_2$ (after perhaps applying complex conjugation). Hence we are reduced to showing $k_1 \cong k_2$.

To that end, let $R_1$ and $R_2$ be the set of finite places of
$k_1$ and $k_2$ that ramify $A_1$ and $A_2$, $R_j'$ denote those
$L$--places lying below those places in $R_j$ for $j=1,2$, and let
$V$ be the set of all finite places of $L$. Let $\omega \in
V\setminus (R_1'\cup R_2')$ be a place that is ramified or inert in
$k_1/L$. Theorem \ref{fuchsiansubgroups} implies that we can
construct a quaternion algebra $B$ over $L$ that is ramified at
$R_1' \cup \omega$ (and possibly other primes) and for which $B
\otimes_L k_1 \cong A_1$. Hence by (\ref{1}), $B$ embeds in $A_2$.
We claim that $\omega$ is ramified or inert in $k_2$. Indeed, if
$\omega=\nu\nu'$ is split in $k_2$, then $(k_2)_\nu \cong L_\omega$
and it follows that $A_2\otimes_{k_2} (k_2)_\nu \cong (B\otimes_L
L_\omega)$ is a division algebra (and similarly for $\nu'$).  Hence
$\nu, \nu' \in R_2$, implying that $\omega\in R_2'$, and this is a
contradiction. Repeating this argument with the roles of $k_1$ and
$k_2$ interchanged we deduce that apart from a finite number of
$L$--primes, an $L$--prime is inert or ramified in $k_1$ if and only
if it is inert of ramified in $k_2$.  It follows then that apart from
a finite number of $L$--primes, the splitting types in $k_1$
and $k_2$ are identical. Proposition \ref{sameplittype} now implies
$k_1=k_2$.\qed\smallskip\smallskip

Note that the proof of Theorem \ref{main3} requires only the knowledge that if a
commensurability class for a surface has a representative in $M_1$,
then it also has a representative in $M_2$ (and conversely). We call
a pair of manifolds $M_1,M_2$ with this property \emph{geometrically
genus commensurable}---compare again with the length spectrum
setting \cite{Re}, \cite{CHLR}, and \cite{PR}.

\begin{theorem}
\label{main3'} Let $M_1={\bf H}^3/\Gamma_1$ and $M_2={\bf
H}^3/\Gamma_2$ be arithmetic hyperbolic 3--manifolds. If $M_1$ and $M_2$ are
geometrically genus commensurable then either they are commensurable or
$\textrm{GS}(M_1)=\emptyset$.
\end{theorem}

%---------------------------------------------------------------
%---------------------------------------------------------------
\section{Aligning surfaces: Theorem \ref{main1}}

In this section we prove Theorem \ref{main1}. This contains most of
the ideas that are needed in the proof of Theorem \ref{main2}. The
remaining complication in the proof of Theorem \ref{main2} will be
controlling multiplicities.

%---------------------------------------------------------------
\subsection{Fuchsian pairs}

The main idea in the proof of Theorem \ref{main1} is given by the
following simple lemma.

\begin{lemma}\label{Basic}
Let $M=\mathbf{H}^3/\Gamma$ be a finite volume hyperbolic
$3$--manifold.\\[\baselineskip]
\noindent (i)~If $\Gamma_1$ and $\Gamma_2$ are (finite index)
subgroups of $\Gamma$ such that
\begin{equation}\label{SurfacePair}
\Gamma_1 \cap F = \Gamma_2 \cap F
\end{equation}
for every maximal Fuchsian subgroup $F$ of $\Gamma$, then the covers
of $M$ corresponding to the subgroups $\Gamma_1$ and $\Gamma_2$ are
geometrically genus equivalent.\\[\baselineskip]
\noindent (ii)~If in addition $G$ is a finite group, $H$ and $K$ are
subgroups of $G$, and $\phi\colon \Gamma \to G$ a surjective
homomorphism such that
\begin{equation}\label{FiniteSurfacePair}
H \cap \phi(F) = K \cap \phi(F)
\end{equation}
for every maximal Fuchsian subgroup $F$ of $\Gamma$, then the covers
of $M$ corresponding to the subgroups $\phi^{-1}(H)$ and
$\phi^{-1}(K)$ are geometrically genus equivalent.
\end{lemma}

\noindent {\bf Proof:}~(ii) follows immediately from (i). That
(i) holds follows, since as noted earlier, the geometric genus set
of $M$ is determined by the maximal Fuchsian subgroups of its associated
Kleinian group.
\qed\\[\baselineskip]

We call a pair of finite groups $H,K<G$ a \emph{Fuchsian pair with
respect to $\phi$} if (\ref{FiniteSurfacePair}) holds.

%---------------------------------------------------------------
\subsection{Geometrical genus equivalence without volume
control}\label{WithoutVolume}

We begin the proof of Theorem \ref{main1} by proving the following
which is nothing more than Theorem \ref{main1} without volume
considerations. The additional volume constraint will be addressed
in Subsection \ref{VolumeControl}.

\begin{theorem}
\label{main1'} Let $M=\mathbf{H}^3/\Gamma$ be an arithmetic
hyperbolic $3$--manifold. Then there exists infinitely many pairs of
nonisometric finite covers $(M_j,N_j)$ of $M$ such that $M_j,N_j$
are geometrically genus equivalent and non-isometric.
\end{theorem}

\noindent \textbf{Proof:}~Without loss of generality we can pass to
a subgroup of finite index so that the invariant trace-field $k$ of
$\Gamma$ coincides with the trace-field. For simplicity,  we will continue to denote this group
by $\Gamma$. Recall that Theorem \ref{fuchsiansubgroups} shows
that the invariant trace-field of any maximal Fuchsian subgroup $F$
of $\Gamma$ coincides with the maximal totally real subfield $L=k
\cap \mathbf{R}$ and $[k:L]=2$.

By the Strong Approximation Theorem (see \cite{Nori87} or
\cite{Weis}) and Lemma \ref{Ceb}, there exists a positive density
set of places $\mathcal{P}$ of $k$ such that for each place $\nu \in
\mathcal{P}$ and place $\omega$ of $L$ with $\nu \mid \omega$,
$\phi_\nu(\Gamma) = \PSL(2,\mathbf{F}_{q_\nu})$, the extension
$\mathbf{F}_{q_\nu}/\mathbf{F}_{q_\omega}$ is proper (that is, it
has degree 2), and the extensions
$\mathbf{F}_{q_\omega}/\mathbf{F}_{p_\nu}$ (the latter being the
prime field) is of some fixed degree $j$. For convenience we shall
omit any dyadic places from $\cal P$.

We will apply Lemma \ref{Basic} in the following way. For each place
$\nu \in \mathcal{P}$, we seek a non-trivial subgroup $C_\nu$ of
$\PSL(2,\mathbf{F}_{q_\nu})$ such that $C_\nu \cap \phi_\nu(F)$ is
trivial for all maximal Fuchsian subgroups $F$ of $\Gamma$. Given
this, the subgroup $C_\nu$ together with the trivial subgroup
$\set{1}$ form a Fuchsian pair with respect to $\phi_\nu$, and so by
Lemma \ref{Basic}, the covers $M_{C_\nu}$ and $M_{1,\nu}$ of $M$
corresponding to the subgroups $\phi_\nu^{-1}(C_\nu)$ and
$\phi_\nu^{-1}(\set{1})$ are geometrically genus equivalent. Furthermore,
note that the manifolds $M_{C_\nu}$ and $M_1$ cannot be isometric since
$M_{1,\nu}$ properly covers $M_{C_\nu}$.

The subgroup $C_\nu$ is provided by the following lemma, the proof
of which is given below.

\begin{lemma}
\label{getsubgroup} There exists a prime $\ell_\nu$ that divides the
order of $\PSL(2,\mathbf{F}_{q_\nu})$ but not the orders of
$\PSL(2,\mathbf{F}_{q_\omega})$ or
$\mathbf{B}(\mathbf{F}_{q_\omega}(\sqrt{b}))$.\end{lemma}

Given the lemma, Cauchy's theorem provides an element $\lambda_\nu$ of order
$\ell_\nu$ in $\PSL(2,\mathbf{F}_{q_\nu})$, and hence a cyclic subgroup $\innp{\lambda_\nu}=C_\nu$
of $\PSL(2,\mathbf{F}_{q_\nu})$ of order $\ell_\nu$. That this
satisfies the condition $C_\nu \cap \phi_\nu(F) = 1$ for any maximal
Fuchsian subgroup follows from Lemma \ref{FuchsianImage} and
elementary group theory.\qed\\[\baselineskip]

\noindent \textbf{Proof of Lemma \ref{getsubgroup}:}~Recall that
\begin{align*}
\abs{\PSL(2,\mathbf{F}_{q_\omega})} &= \frac{q_\omega(q_\omega^2-1)}{2} = \frac{q_\omega(q_\omega+1)(q_\omega-1)}{2}, \\
\abs{\mathbf{B}(\mathbf{F}_{q_\omega}(\sqrt{b}))} &=\frac{q_\omega^2(q_\omega^2-1)}{2} = \frac{q_\omega^2(q_\omega+1)(q_\omega-1)}{2} \\
\abs{\PSL(2,\mathbf{F}_{q_\nu})} &=\frac{q_\omega^2(q_\omega^2+1)(q_\omega^2-1)}{2} = \frac{q_\nu(q_\nu+1)(q_\nu-1)}{2}.
\end{align*}
Since the numbers $q_\nu+1,q_\omega$ are relatively prime and the
numbers $q_\nu-1,q_\nu+1$ share only 2 as a common prime divisor,
Lemma \ref{getsubgroup} can be established by finding an odd prime
divisor of $q_\nu+1$. However, this is elementary as $q_\nu+1 =
2\alpha_\nu$ where $\alpha_\nu \geq 5$ is odd. Thus \emph{any} prime
divisor $\ell_\nu$ of $\alpha_\nu$ will suffice for the
lemma.\qed\\[\baselineskip]

\noindent \textbf{Remarks:}\smallskip\smallskip

\noindent (1)~There is an obvious similarity between
the statement of Lemma \ref{Basic} and a basic lemma used in the
methods of Sunada \cite{Sun} and \cite{LMNR06}. Thus our method can
be viewed as a generalization of Sunada's method to the setting of geometric
genus spectra. However, unlike the case of Sunada our method is not
purely algebraic. The reason is this: in Sunada's method, $\pi_1(S^1)$ is
cyclic, so under any homomorphism $\phi$ of $\Gamma$, the image of
$\pi_1(S^1)$ is cyclic, and so $\phi$ plays a minimal role while
the role of $G$ is paramount. In contrast, in applying Lemma \ref{Basic},
one must have knowledge of the image of Fuchsian subgroups of
$\Gamma$ under $\phi$, images which could be extremely complicated. \\[\baselineskip]
\noindent (2)~Although Lemma \ref{FuchsianImage} is stated only for totally
geodesic surfaces, one can formulate a statement for surface
subgroups of Kleinian groups more generally. However, there seems to
be no way of controlling where the image of a general surface group
maps under any homomorphism onto a finite group. Indeed, when the trace-field of the
surface group coincides with that of the Kleinian group (which is the typical situation when
the surface subgroup is not Fuchsian), the Strong Approximation
Theorem ensures that the restriction of $\phi_\nu$ to the surface
group is surjective for all but finitely many places.

The terminology \emph{geometric} genus spectrum is meant to
emphasize that the free homotopy classes of surfaces in the
3--manifold under consideration are of a geometric origin. The
\emph{genus spectrum} of $M$ is the set of pair
$((g,n),\ell_{g,n}(M))$ of free homotopy classes of
$\pi_1$--injective, properly immersed surfaces of $M$ that contain a
representative of topological type $\Sigma_{g,n}$. By \cite{T}
Corollary 8.8.6, if $M$ is a closed hyperbolic 3--manifold then
again $\ell_{g,n}(M)$ is finite, and it would seem to be an interesting
problem to see how much of the topology and geometry of $M$ is
determined by this topological spectrum.

%---------------------------------------------------------------
\subsection{Ensuring volume growth: Proof of Theorem
\ref{main1}}\label{VolumeControl}

We now show how to extend the method of proof of Theorem \ref{main1'}  to obtain
Theorem \ref{main1}. \smallskip\smallskip

\noindent \textbf{Proof of Theorem \ref{main1}:}~ The proof of Lemma
\ref{getsubgroup} shows that any odd prime divisor $\ell_\nu$ of
$q_\nu+1$ produces a pair of geometrically genus equivalent covers
with volume ratio $\ell_\nu$. Ranging over $\mathcal{P}$, if the set
of $\ell_\nu$ is unbounded, then the result follows.
Otherwise, there is a finite list of primes $\ell_1,\dots,\ell_r$
such that (recall $q_\nu+1$ is never zero modulo 4)
\[ q_\nu+1 = 2\prod_{i=1}^r \ell_i^{\alpha_{i,\nu}}. \]
Note that $q_\nu$ and the exponents $\alpha_{i,\nu}$ depend on the
place $\nu \in \mathcal{P}$ but the primes $\ell_1,\dots,\ell_r$, by
assumption, do not. In particular, since the left hand side is
unbounded as we range over $\mathcal{P}$, one of the exponents
$\alpha_{i_0,\nu}$ must be unbounded as we range over $\mathcal{P}$.
According to Lemma \ref{OrderLemma}, since $\ell_i^{\alpha_{i_0,\nu}}$
divides $q_\nu+1$, there exists an element $\lambda_\nu$ of order
$\ell_i^{\alpha_{i_0,\nu}}$ in $\PSL(2,\mathbf{F}_{q_\nu})$. Arguing as before, we see that the
subgroup $C_\nu$ generated by $\lambda_\nu$ and the trivial subgroup form a
Fuchsian pair for $\phi_\nu$. The resulting geometrically genus
equivalent covers of $M$ have a volume ratio of
$\ell_i^{\alpha_{i_0,\nu}}$, which by choice of
$\alpha_{i_0,\nu}$, is unbounded as we range over
$\mathcal{P}$.\qed\\[\baselineskip]

\noindent \textbf{Remark:}~The manifolds $M_{C_\nu}$ and $M_{1,\nu}$
in the proof of Theorem \ref{main1'} are \emph{not} geometrically isogenus
when $\textrm{GS}(M) \ne \emptyset$.
To see this, note that the maximal orientable totally geodesic
surfaces up to free homotopy in a hyperbolic 3--manifold
$M=\mathbf{H}^3/\Gamma$ are in bijection with the
$\Gamma$--conjugacy classes of maximal Fuchsian subgroups of
$\Gamma$. Thus given a finite cover $N \to M$, the maximal
orientable totally geodesic surfaces up to free homotopy in $N$ are
parameterized by the $\pi_1(N)$--conjugacy classes of the Fuchsian
subgroups $\pi_1(N) \cap F$, where $F$ is a maximal Fuchsian
subgroup of $\Gamma$.

Applying this discussion to the manifolds $M_{C_\nu}$ and
$M_{1,\nu}$, we know that $\pi_1(M_{C_\nu})$ and  $\pi_1(M_{1,\nu})$
satisfy (\ref{SurfacePair}). Since $M_{1,\nu} \to M_{C_\nu}$ is a
cyclic cover of degree $\ell_\nu$ and (\ref{SurfacePair}) holds,
each $\pi_1(M_{C_\nu})$--conjugacy class $\pi_1(M_{C_\nu}) \cap F$
produces $\ell_\nu$ distinct $\pi_1(M_{1,\nu})$--conjugacy classes
in $\pi_1(M_{1,\nu})$. This yields the relationship:
\begin{equation}\label{GenusFormula}
{\cal GS}(M_{\nu,1}) = \set{ (X,\ell_X(M_{\nu,1})} =
\{(X,\ell_\nu\cdot \ell_X(M_{\nu,\ell_\nu}))\}
\end{equation}
In particular, \emph{none} of the non-zero multiplicities are the same in the
geometric genus spectra. However, the uniform nature of this failure
will provide us a handle for matching up multiplicities for
other pairs of covers.

%---------------------------------------------------------------
%---------------------------------------------------------------
%---------------------------------------------------------------
\section{Arranging multiplicities: Proof of Theorem \ref{main2}}

In this section we prove Theorem \ref{main2}. As discussed above,
the proof of Theorem \ref{main1'} produces infinitely many pairs of
finite covers $M_{C_\nu}$ and $M_{1,\nu}$ that are geometrically
genus equivalent with a very precise relationship between their
geometric genus spectra given by (\ref{GenusFormula}). We now show
how to exploit (\ref{GenusFormula}) to produce infinitely many pairs
of geometrically isogenus manifolds by using product homomorphisms
$\phi_{\nu_1} \times \phi_{\nu_2}$. More precisely, Theorem
\ref{main2} follows from (using the notation from the proof of Theorem
\ref{main1'}):

\begin{proposition}
\label{commonell} There exists an infinite subset of places $\nu\in
{\cal P}$ for which a fixed prime $\ell$ can be taken for
$\ell_\nu$.\end{proposition}

Deferring the proof of this we complete the proof of Theorem \ref{main2}.\\[\baselineskip]

\noindent \textbf{Proof of Theorem \ref{main2}:} \noindent Let
$\mathcal{P}'\subset{\cal P}$ be the subset produced by Proposition
\ref{commonell} for the prime $\ell$. For any pair of places $\nu_1,\nu_2 \in
\mathcal{P}'$, we have a pair of surjective homomorphisms
$\phi_{\nu_j}\colon \Gamma \to \PSL(2,\mathbf{F}_{q_{\nu_j}})$ such
that $\ell$ divides the order of $\PSL(2,\mathbf{F}_{q_{\nu_j}})$
but not the orders of $\PSL(2,\mathbf{F}_{q_{\omega_j}})$ or
$\mathbf{B}(\mathbf{F}_{q_{\omega_j}}(\sqrt{b}))$. Let $C_{\nu_j}$
be a cyclic subgroup of $\PSL(2,\mathbf{F}_{q_{\nu_j}})$ of order
$\ell$. Setting
\[ \phi_{1,2}\colon \Gamma \to \PSL(2,\mathbf{F}_{q_{\nu_1}}) \times \PSL(2,\mathbf{F}_{q_{\nu_2}}) \]
to be the product homomorphism $\phi_{\nu_1} \times \phi_{\nu_2}$
(which is also surjective), we define subgroups $H_1= C_{\nu_1}
\times \set{1}$ and $H_2=\set{1} \times C_{\nu_2}$. By construction
of $\ell$, it follows that $H_1,H_2$ satisfy (\ref{FiniteSurfacePair}),
and thus form a Fuchsian pair for $\phi_{1,2}$. Hence, by Lemma
\ref{Basic}(ii), the covers
$M_{\nu_1,\nu_2,H_1},M_{\nu_1,\nu_2,H_2},M_{\nu_1,\nu_2,\set{1}}$ associated
to the subgroups
$\phi_{1,2}^{-1}(H_1),\phi_{1,2}^{-1}(H_2),\phi_{1,2}^{-1}(\set{1})$ are
pairwise geometrically genus equivalent. Indeed, by
(\ref{GenusFormula})
\[ {\cal GS}(M_{\nu_1,\nu_2,\set{1}}) = \{(X,\ell\cdot \ell_X(M_{\nu_1,\nu_2,H_j}))\} \]
for $j=1,2$. In particular,
\[ \ell\cdot \ell_X(M_{\nu_1,\nu_2,H_1}) = \ell \cdot \ell_X(M_{\nu_1,\nu_2,H_2}) \]
and thus
\[ \ell_X(M_{\nu_1,\nu_2,H_1}) = \ell_X(M_{\nu_1,\nu_2,H_2}). \]

It remains to prove that $M_{\nu_1,\nu_2,H_1},M_{\nu_1,\nu_2,H_2}$ are
not isometric. If this were not the case, then $\pi_1(M_{\nu_1,\nu_2,H_1})$
and $\pi_1(M_{\nu_1,\nu_2,H_2})$ would be conjugate in
$\textrm{Isom}(\mathbf{H}^3)$. We claim this is impossible. To that end, set
$\lambda_{\nu_1} \in C_{\nu_1}$ to be a generator and $\gamma_{\nu_1} \in
\pi_1(M_{\nu_1,\nu_2,H_1})$ such that $\phi_1(\gamma_{\nu_1}) = \lambda_{\nu_1}$.
Recall that the homomorphisms $\phi_\nu$ arose from the reduction of
$R_\nu$ be the unique maximal ideal $\pi R_\nu$ and so induce
homomorphisms $\phi_\nu\colon R_\nu[t] \to \mathbf{F}_{q_\nu}[t]$.
It is a simple matter that for $\eta \in \PSL(2,R_\nu)$, we have the
relationship $c_{\phi_\nu(\eta)}(t) = \phi_\nu(c_\eta(t))$, where
$c_\theta(t) \in R_\nu[t]$ denotes the characteristic polynomial of
$\theta$.

By definition, every element $\gamma$ of $\pi_1(M_{\nu_1,\nu_2,H_2})$
maps trivially under $\phi_{\nu_1}$ and thus
$\phi_{\nu_1}(c_{\gamma}(t)) = (t-1)^2$. As conjugation by any
$\tau$ in $\textrm{Isom}(\mathbf{H}^3)$ preserves determinant and at
worst changes trace by complex conjugation, any
$\textrm{Isom}(\mathbf{H}^3)$--conjugate of $\gamma$ also has this
property. In particular, if $\gamma_{\nu_1}$ is
$\textrm{Isom}(\mathbf{H}^3)$--conjugate into
$\pi_1(M_{\nu_1,\nu_2,H_2})$, it would have to be that $c_{\lambda_{\nu_1}}(t) =
(t-1)^2$. However, this is impossible since $\lambda_{\nu_1}$ is a
semisimple element of odd prime order $\ell$.

To produce infinitely many pairs, we can simply vary $\nu_1,\nu_2$ over $\mathcal{P}'\times \mathcal{P}' \setminus \Delta$
where $\Delta$ denotes the diagonal.
\qed\\[\baselineskip]

\noindent \textbf{Proof of Proposition \ref{commonell}:} In the proof
of the proposition, we will keep with the notation used in the proof
of Theorem \ref{main1'}. For the reader's convenience, we will briefly
recall some of the notation here. Recall that we restricted ourselves to a positive density subset
${\cal P}$ of places of $k$ in order to control
$[\mathbf{F}_{q_\nu}:\mathbf{F}_{q_\omega}]$ and $[\mathbf{F}_{q_\omega}: \mathbf{F}_{p_\nu}]$
and will continue to do so here. In addition, $q_\nu,q_\omega,p_\nu$ will denote the orders of the field
$\mathbf{F}_{q_\nu},\mathbf{F}_{q_\omega}$, and $\mathbf{F}_{p_\nu}$, and  $\textrm{P}$ will
denote the set of all integral primes.

For each $\nu\in {\cal P}$, we saw from the proof of Theorem \ref{main1'} that
$\ell_\nu$ can be any odd prime divisor of $q_\nu+1$. What is required here is to find an odd prime $\ell$ and an infinite
subset ${\cal P}' \subset {\cal P}$ such that $\ell \mid q_\nu+1$ for all $\nu \in {\cal P}'$. In fact,
we will find an odd prime $\ell$ and a positive density (and hence infinite) subset ${\cal P}'$ for which this holds. The remainder of this
proof is devoted to this task.

Recall that $q_\nu= q_\omega^2=p_\nu^{2j}$ where $j$
is fixed, and so $q_\nu+1 = p_\nu^{2j}+1$. It is elementary that for any
(odd) prime $\ell$, if $\ell$ divides $p_\nu^{2j}+1$, then $p_\nu^{2j}+1 \equiv 0 \mod \ell$. Equivalently,
setting $F_j(t) = t^{2j}+1$, the previous statement is simply that  $\overline{p_\nu}$ is a zero of $F_j(t)$ in
$\mathbf{F}_\ell$ (where $\overline{p_\nu}$ is the modulo $\ell$ residue class of $p_\nu$ in $\mathbf{F}_\ell$).

We now show how to
use $F_j(t)$ and the Cebotarev density theorem to find the required $\ell$ and ${\cal P}'$. To that end, let $\mathcal{L}_j$ be the set
of odd primes $\ell$ such that $F_j(t)$ has a root in
$\mathbf{F}_\ell$. For each $\ell \in \mathcal{L}_j$,
define the set
\[ \textrm{P}^j_\ell = \set{m \in \mathbf{N}~:~ F_j(\overline{m}) = 0 \mod \ell} \cap \textrm{P}. \]
By the Cebotarev Density Theorem,
\begin{equation}\label{PositiveDensity}
\textrm{Density}(\textrm{P}^j_\ell) \geq \frac{1}{\ell-1}
\end{equation}
In addition, since for every prime $p$, there exists an odd prime divisor of $p^{2j}+1$
(recall that $p^{2j}+1$ is never zero modulo 4),
\begin{equation}\label{Density1}
\bigcup_{\ell \in \mathcal{L}_j} \textrm{P}^j_\ell = \textrm{P}.
\end{equation}
After ordering $\mathcal{L}_j = \set{\ell_1<\ell_2<\ell_3<\dots}$,
for each $i \geq 1$, we define the set
\[ A_i = \bigcup_{m=1}^i \textrm{P}^j_{\ell_m}. \]
The sequence of densities $\set{\textrm{Density}(A_i)}_i$ for the
sets $A_i$ is positive by (\ref{PositiveDensity}), strictly
increasing (since $\textrm{P}^j_{\ell_i} \cap \textrm{P}^j_{\ell_{i'}}$ has density strictly smaller than
$\textrm{P}^j_{\ell_i}$ and $\textrm{P}^j_{\ell_{i'}}$ for all distinct $i,i'$), and has a least upper bound of $1$
by (\ref{Density1}). This in tandem with the fact that the set
\[ \textrm{P}_{\mathcal{P}} = \set{p_\nu~:~ \nu\in {\cal P}} \]
has positive density, implies that there exists an $i_0$ such that
$A_{i_0} \cap \textrm{P}_{\mathcal{P}}$ has positive density. By the Pigeonhole Principle,
there exists $1 \leq i_1 \leq i_0$ such that
\[ \textrm{P}' = \textrm{P}^j_{\ell_{i_1}} \cap \textrm{P}_{\mathcal{P}} \]
has positive density. In particular, for each $p \in \textrm{P}'$, the
prime $\ell_{i_1}$ divides $p^{2j}+1$. Thus, $\ell_{i_1}$ and the set
\[ \mathcal{P}' = \set{\nu \in \mathcal{P}~:~ p_\nu \in \textrm{P}'} \]
is the required pair needed to finish the proof.\qed

%---------------------------------------------------------------
%---------------------------------------------------------------
\section{Examples}

In this section, we discuss our constructions in the context of some
examples.

%---------------------------------------------------------------
\subsection{} 

Take $k=\mathbf{Q}(i)$ and $L=\mathbf{Q}$.
Since any prime $p \equiv 3 \mod 4$ is inert in $k$, we see that
$\mathbf{F}_\nu \ne \mathbf{F}_p$ for the unique place $\nu$ dividing
$p$. Consider the primes $3$ and
$7$. Note that $\abs{\PSL(2,\mathbf{F}_9)} = 360$ and
$\abs{\PSL(2,\mathbf{F}_{49})} = 58800$ are both divisible by $5$.
On the other hand $\abs{\PSL(2,\mathbf{F}_3)} = 12$,
$\abs{\PSL(2,\mathbf{F}_7)} = 168$, $\abs{\mathbf{B}(\mathbf{F}_9)}
= 36$ and $\abs{\mathbf{B}(\mathbf{F}_{49})} = 1176$, all fail to have $5$
as a divisor.

With this, let $A/k$ be any quaternion algebra which satisfies the hypothesis
of Theorem \ref{fuchsiansubgroups} (and so contains Fuchsian subgroups)
and for which the places
of $k$ above $7$ are not contained in $\textrm{Ram}_f A$. Note that by
definition of $A$, the place above $3$ is not contained in $\textrm{Ram}_f A$.
Let $\cal O$ be an (maximal) order of $A$ and $\Gamma$ the image in
$\PSL(2,{\bf C})$ of ${\cal O}^1$.  In the notation of the proof of Theorem \ref{main1'},
we can take $\nu$ to be the place above $3$, and $\ell_\nu=5$ (so that
$C_\nu$ is a cyclic group of order $5$). The proof of Theorem
\ref{main1'} now produces geometrically genus equivalent covers of
degree $72$ and $360$, respectively.

To produce isogenus covers, we can take the places over $3$ and
$7$, and use $\ell=5$. The proof of Theorem \ref{main2} now produces geometrically
isogenus covers of degree $4,233,600$.

We remark that the reduction
homomorphisms do indeed surject $\Gamma$ onto the stated finite groups (without recourse
to Strong Approximation). The reason for this
is that since $\cal O$ is maximal, it will be dense in its localizations. Since $\nu$ is unramified,
these localizations of $\cal O$ are isomorphic to $M(2,{\bf Z}[i]_\nu)$, and thus $\Gamma$
is dense in the groups
$\PSL(2,{\bf Z}[i]_\nu)$ for the stated places $\nu$.

Note that there are infinitely many distinct commensurability
classes of such examples that satisfy this condition on
$\textrm{Ram}_f A$. Indeed we can arrange for $\Gamma$ to be
torsion-free.

%---------------------------------------------------------------
\subsection{} 

Using the results of \cite{JM} and
\cite{MR1} one can in principle obtain information on the
topological types and multiplicities arising in the case of the
groups $\PSL(2,O_d)$. These groups contain elements of finite order,
but the arguments still apply. We sketch some of this for the case
Picard group $\PSL(2,O_1)$.

Maximal arithmetic Fuchsian groups are parameterized by a positive
integer $D$, called the discriminant.  In \cite{MR1} it is shown
that the number of $\PSL(2,O_1)$--conjugacy classes of maximal
arithmetic Fuchsian subgroups of discriminant $D$ is 1, 2, or 3
dependent on whether $D$ is congruent to 0 or 3 mod 4, 2 mod 4 or 1
mod 4. Indeed, one can get more detailed information on these
arithmetic Fuchsian groups.

For all discriminants $D$, one conjugacy class is represented by the
arithmetic Fuchsian group obtained as follows. Let $B_D$ be the
indefinite quaternion algebra over $\bf Q$ with Hilbert Symbol
$\biggl({{-1,D}\over {\bf Q}}\biggr)$ and ${\cal O}_D = {\bf
Z}[1,i,j,ij]$. Then prescribing a particular representation
$\rho:B_D\to M(2,{\bf C})$, the image of ${\cal O}_D^1$ determines
an arithmetic Fuchsian subgroup which we denote by $F_D$. When $D$
is congruent to $0$ or $3$ modulo $4$ this is the only such group
(up to $\PSL(2,O_1)$--conjugacy).

In the case when $D$ is congruent to 1 modulo $4$, the two further
conjugacy classes are represented by group $G_{D,1}$ and $G_{D,2}$.
As remarked in \cite{MR1}, the subgroups $G_{D,1}$ and $G_{D,2}$ are
conjugate in $\PGL(2,O_1)$. It is shown in \cite{MR1} that the group
$G_{D,1}$ is the image in $\PSL(2,O_1)$ of the unit group of an
order ${\cal M}_D$ of $B_D$ and $[{\cal M}_D^1:{\cal O}_D^1]=3$.

Finally for the case when $D$ is congruent to $2$ modulo 4, there is a
further group $H_D$ that is described as the image of the unit group
of another order ${\cal N}_D$ of $B_D$ and in this case $[{\cal
N}_D^1:{\cal O}_D^1] = 2$ or $6$ dependent on whether $D$ is
congruent to $2$ or $6$ modulo 8.

From this the co-volumes of these maximal arithmetic Fuchsian groups
can be computed using a result of Humbert (see \cite{MR1} \S 6). In
addition, whether the group is non-cocompact is also decidable from
$D$ (see Lemma 8 of \cite{MR1}) and the numbers of conjugacy classes
of elements of finite order (and the orders) are also computable
(see \cite{MR1} \S 7).

%---------------------------------------------------------------
%---------------------------------------------------------------
\section{The non-arithmetic case}\label{Non-Arithmetic}

In this section, we outline how Theorems \ref{main2} and
\ref{main1} can be extended to any finite volume hyperbolic
3--manifold.

Let $M=\mathbf{H}^3/\Gamma$ be a complete, finite volume hyperbolic
3--manifold. Set $k$ to be the invariant trace field of $M$ and $F$
to be the maximal real subfield. For each place $\nu$ of $k$ and
$\omega$ of $F$ with $\nu \mid \omega$, let
$\mathbf{F}_{q_\nu},\mathbf{F}_{q_\omega}$ be the associated residue
fields and $t_\nu=[\mathbf{F}_{q_\nu}:\mathbf{F}_{q_\omega}]$. According
to the Strong Approximation Theorem and Lemma \ref{Ceb}, we can
arrange it so that $t=t_\nu>1$ is constant, $[\mathbf{F}_{q_\omega}:\mathbf{F}_{p_\nu}]=j$ is constant, and $\phi_\nu(\Gamma) =
\PSL(2,\mathbf{F}_{q_\nu})$ for a positive density subset $\mathcal{P}$ of places $\nu$ of $k$.

The proof of Theorem \ref{main1'} for general $M$ is done with
precisely the same method. That Lemma  \ref{getsubgroup} holds in this setting follows immediately
from Zsigmondy's theorem (\cite{Zs}).  However, to prove the analogue of Proposition
\ref{commonell} in this setting requires more than what Zsigmondy's theorem obviously gives, and with that
in mind we therefore give more discussion. First, note that one can simplify things somewhat. For each $\nu \in \mathcal{P}$,
take a subfield $\mathbf{F}_{q_\omega} < \mathbf{F}_\nu < \mathbf{F}_{q_\nu}$ such that the
degree $[\mathbf{F}_\nu:\mathbf{F}_{q_\omega}] = p'$ is a fixed prime. Notice that it
certainly suffices to find a prime divisor $\ell_\nu$ of $|\PSL(2,\mathbf{F}_\nu)|$ that does not divide
$|\PSL(2,\mathbf{F}_{q_\omega})|$ or $|\mathbf{B}(\mathbf{F}_{q_\omega}(\sqrt{b}))|$ as this would
also divide the order of the larger group $\PSL(2,\mathbf{F}_{q_\nu})$.
For $p'=2$, this is nothing more than Lemma \ref{getsubgroup}. However, when $p'$ is odd, we need to amend
the  proof of Proposition \ref{commonell} as follows.
In this case,  we can write
\[ |\PSL(2,\mathbf{F}_\nu)| = \frac{q_\omega^{p'} (q_\omega+1)(q_\omega-1) r_{-,\nu} r_{+,\nu}}{2} \]
where the factors $r_{\pm,\nu}$ are given by:
\begin{equation}\label{OddPolys}
r_{+,\nu} = \frac{q_\omega^{p'}+1}{q_\omega+1}, \quad r_{-,\nu} =
\frac{q_\omega^{p'}-1}{q_\omega-1},
\end{equation}
and have the following divisibility properties:
\begin{equation}\label{Div}
(q_\omega\pm 1, r_{\mp,\nu}) = 1,\quad (r_{\pm,\nu},q_\omega)=1, \quad (r_{+,\nu},r_{-,\nu})=1, \quad (r_{\pm, \nu},q_{\omega}\pm 1) \mid p'.
\end{equation}
Assuming the validity of (\ref{Div}), the existence of $\ell_\nu$ is easy. As either $r_{+,\nu}$ or $r_{-,\nu}$ must be relatively prime to $p'$ and both have at least one odd prime divisor, there exists an odd prime divisor of either $r_{+,\nu}$ or $r_{-,\nu}$ that does not divide $|\PSL(2,q_{\omega})|$ and $|\mathbf{B}(\mathbf{F}_{q_\omega}(\sqrt{b}))|$.

The proof of (\ref{Div}) is straightforward. In fact, the first three assertions are trivial and only the proof of the last assertion requires comment. To begin, an elementary calculation shows that
\[ r_{-,\nu} = \sum_{j=0}^{p'-1} q_\omega^j, \quad r_{+,\nu} = \sum_{j=1}^{p'} (-1)^{p'-j}q_\omega^{j-1}. \]
Using this, we will work out the case of $(r_{-,\nu},q_\omega-1)$. Let $\ell$ be a common divisor of $r_{-,\nu}$ and $q_\omega-1$. Then $\ell$ must divide
\[ \pr{r_{-,\nu} + (q_\omega-1)} = \pr{\sum_{j=0}^{p'-1} q_\omega^j} + q_\omega - 1 = 2q_\omega + \sum_{j=2}^{p'-1} q_\omega^j. \]
From this, we see that that $\ell$ must divide
\[ \pr{\pr{ 2q_\omega + \sum_{j=2}^{p'-1} q_\omega^j} + 2q_\omega(q_\omega-1)} = 3q_\omega^2 + \sum_{j=3}^{p'-1}q_\omega^j. \]
Using this procedure, we can now inductively deduce:
\begin{equation}\label{MinusInduct}
\ell \mid  \pr{ mq_\omega^{m-1} + \sum_{j=m}^{p'-1} q_\omega^j}
\end{equation}
for any integer $2 \leq m \leq p'-1$. Taking $m=p'-1$ in (\ref{MinusInduct}) yields that $\ell$ must divide
\[ (1-p')q_\omega^{p'-2} + q_\omega^{p'-1} = q_\omega^{p'-2}(q_\omega + p'-1) \]
and thus $\ell \mid (q_\omega+p'-1)$. This in tandem with the fact that $\ell \mid (q_\omega-1)$ implies that $\ell \mid p'$ as desired. For the case of $\ell \mid (r_{+,\nu},q_\omega+1)$, the corresponding inductive statement becomes:
\begin{equation}\label{PositiveInduct}
\ell \mid \pr{ (-1)^{p'-m} m q_\omega^{m-1} + \sum_{j=m+1}^{p'} (-1)^{p'-j} q_\omega^{j-1}}
\end{equation}
for all $2 \leq m \leq p'-1$. Taking $m=p'-1$ in (\ref{PositiveInduct}) implies that $\ell$ divides
\[ (-p'+1)q_\omega^{p'-2} + q_\omega^{p'-1} = q_\omega(1-p' + q_\omega). \]
As before, this and $\ell \mid q_\omega+1$ implies that $\ell \mid p'$.

To prove the extension of Theorem \ref{main1} in this case, what is important here is that the prime $\ell_\nu$
obtained from the extension of Lemma \ref{getsubgroup} occurs as a
prime divisor of one of $r_{+,\nu},r_{-,\nu}$. Taking a positive
density subset of places $\nu$ where the $\ell_\nu$ arises as a
factor of a fixed $r_{\star,\nu}$, the argument is identical.

For the extension of Theorem \ref{main2}, we replace the polynomial $x^{2j}+1$ by one
associated to either $r_{+,\nu}$ or $r_{-,\nu}$  obtained from (\ref{OddPolys}) (at least one occurs
with positive density). With this, the
argument is identical to that given in the proof of Proposition \ref{commonell}.

%---------------------------------------------------------------
%---------------------------------------------------------------
\section{Final remarks and questions}\label{Final}

In this final section, we collect some questions that naturally arise from this work.

%---------------------------------------------------------------
\subsection{Geometric genus spectrum invariants}

All of the examples of geometrically isogenus manifolds produced in
Theorem \ref{main2} have equal volume. As this is a well known
spectral invariant, this prompts the question:\smallskip\smallskip

\noindent \textbf{Question.} \emph{Do geometrically isogenus
(arithmetic) hyperbolic 3--manifolds with non-trivial geometric genus spectra
always have equal volume?}

%---------------------------------------------------------------
\subsection{A criteria for arithmeticity}

As Theorem \ref{main1} shows, any pair of arithmetic hyperbolic
3--manifolds with identical, nontrivial geometric genus sets are
commensurable. As already discussed, we do not know of a non-arithmetic
example for which geometric genus set is
infinite.\smallskip\smallskip

\noindent \textbf{Question.} \emph{If the geometric genus set of a
complete, finite volume hyperbolic 3--manifold $M$ is infinite, is
$M$ arithmetic?}

%---------------------------------------------------------------
\subsection{Commensurability in general}

As we expect non-arithmetic complete, finite volume hyperbolic
3--manifolds to typically have finite geometric genus set, a
generalization of Theorem \ref{main3} to non-arithmetic manifolds
seems unlikely. In particular, we expect an affirmative answer to
the following question.\smallskip\smallskip

\noindent \textbf{Question.} \emph{Do there exist incommensurable
non-arithmetic 3--manifolds with equal, nontrivial geometric genus
sets? What about equal nontrivial geometric genus
spectra?}\smallskip\smallskip

A class of examples which are good candidates for an affirmative
answer to these questions are certain non-arithmetic hyperbolic
twist knot complements. Any hyperbolic twist knot complement
contains an immersed twice punctured disk which is always totally
geodesic with a unique hyperbolic structure (arising from the level
2 principal congruence subgroup of $\PSL(2,{\bf Z})$). Now using
\cite{HS} we can arrange for infinitely many of these twist knot
complements to have invariant trace-fields of odd prime degree.  It
now follows from Theorem 5.3.8 of \cite{MR} that these knot
complements contain no closed totally geodesic surfaces, and
furthermore the proof shows that any totally geodesic surface in the
knot complement covers ${\bf H}^2/\PSL(2,{\bf Z})$.

At present we cannot rule out the existence of other surfaces but we
believe that the totally geodesic twice punctured disk (and its
covers) are the only totally geodesic surfaces in these knot
complements.

%---------------------------------------------------------------
%---------------------------------------------------------------

%---------------------------------------------------------------
%---------------------------------------------------------------

\noindent Department of Mathematics \\
University of Chicago \\
Chicago, IL 60637 \\
email: {\tt dmcreyn@math.uchicago.edu} \\

\noindent Department of Mathematics\\
University of Texas\\
Austin, TX 78712\\
email: {\tt areid@math.utexas.edu}

%---------------------------------------------------------------
%---------------------------------------------------------------

\end{document}